\documentclass[11pt]{article}
\usepackage{amsmath,amsthm,latexsym,amssymb}
\usepackage{amsfonts}
\usepackage{amsbsy}

\def\be{\begin{eqnarray}}
\def\ee{\end{eqnarray}}

\usepackage{amsmath}
\usepackage{amsfonts}
\usepackage{amsbsy}

\title{\textbf{Measurability of optimal transportation and
convergence rate for Landau type interacting  particle systems}}
\author{Joaquin Fontbona\thanks{DIM-CMM, Universidad de Chile, Casilla 170-3, Correo 3, Santiago-Chile,
e-mail:fontbona@dim.uchile.cl. Supported by Fondecyt Proyect
1040689, ECOS-Conicyt C05E02, Millennium Nucleus Information and
Randomness ICM P04-069-F and FONDAP Applied Mathematics},
H\'el\`ene Gu\'erin
\thanks{IRMAR, Universit\'e Rennes 1, Campus de Beaulieu, 35042 Rennes-France, e-mail:helene.guerin@univ-rennes1.fr.
Supported by  ECOS-Conicyt C05E02 and Millennium Nucleus
Information and Randomness ICM P04-069-F}, Sylvie
M\'el\'eard\thanks{CMAP, Ecole Polytechnique, CNRS, route de
Saclay, 91128 Palaiseau Cedex-France e-mail:
sylvie.meleard@polytechnique.edu.  Supported by ECOS-Conicyt
C05E02 and Millennium Nucleus Information and Randomness ICM
P04-069-F}}

\newtheorem{theo}{\bf Theorem}[section]
\newtheorem{defi}[theo]{\bf Definition}
\newtheorem{lem}[theo]{\bf Lemma}
\newtheorem{prop}[theo]{\bf Proposition}
\newtheorem{coro}[theo]{\bf Corollary}
\newtheorem{rem}[theo]{\bf Remark}

\def\RR{\mathbb{R}}
\def\NN{\mathbb{N}}
\def \EE{\mathbb{E}}
\def \PP{\mathbb{P}}

\def \1{\mathbf{1}}

\def\be{\begin{eqnarray}}
\def\ee{\end{eqnarray}}

\date{}
\renewcommand{\Box}{\hfill\mbox{\fbox{\rule{0mm}{1.5mm}}}}

\begin{document}

\pagenumbering{arabic}
\maketitle
 \pagestyle{plain} \thispagestyle{empty}

\begin{abstract}
In this paper, we consider  nonlinear  diffusion processes driven
by space-time white noises, which have  an interpretation in terms
of partial differential equations. For a  specific choice of
coefficients, they correspond to  the Landau equation arising in
kinetic theory. A particular feature is that the diffusion matrix
of this process is a linear function the law of the process, and
not a quadratic one, as in the McKean-Vlasov model. The main goal
of the paper is to construct an easily simulable diffusive
interacting particle system, converging towards this nonlinear
process and to obtain an explicit pathwise rate.
 This  requires to find a significant coupling between finitely many Brownian
motions and the infinite dimensional white noise process. The key
idea will be  to construct the right Brownian motions by pushing
forward the white noise processes, through the Brenier map
realizing the optimal transport between the law of the nonlinear
process, and the empirical measure of independent copies of it. A
striking problem then is to establish the joint measurability of
this optimal transport map with respect to the space variable  and
the parameters (time and randomness) making the marginals vary. We
shall prove a general  measurability result for the  mass
transportation problem in terms of the support of the transfert
plans,  in the sense of set-valued mappings. This will allow us to
construct the coupling and to obtain explicit convergence rates.
\end{abstract}



\bigskip

\section{Introduction and main statements}
Consider the nonlinear diffusion processes in $\RR^d$ of the
following
 type:
\begin{equation}\label{eq:nonlinproc}
X_t=X_0+\int_0^t \int_{\RR^d} \sigma(X_s-y) W_P(dy,ds)+\int_0^t
\int_{\RR^d} b(X_s-y) P_s(dy)ds
\end{equation}
where $P_t$ is the law  of  $X_t$, and $\ W_P\ $ is  a $\RR^d$
valued space-time white noise on $[0,T]\times \RR^d$ with
independent coordinates, each of which having covariance measure
$P_t(dy)\otimes dt$. \smallskip

 The nonlinear process
(\ref{eq:nonlinproc}) was introduced by Funaki \cite{Funaki:84},
who obtained existence and uniqueness results for Lipschitz
coefficients $\sigma:\RR^d\to \RR^{d\otimes d}$ and $b:\RR^d\to
\RR^d$, see also  Guerin \cite{Guerin:03} for a different
approach. It has an important interpretation in terms of
 partial differential equations issued from kinetic theory. More precisely, for a specific choice of coefficients
$\sigma$ and $b$, the laws $(P_t)_t$ are a weak solution of the
spatially homogeneous
 Landau (also called Fokker-Planck-Landau)
 equations for Maxwell potential:
 \begin{equation}\label{eq:Landaueq}
\frac{\partial f}{\partial t}\left( t,v\right) =\frac{1}{2}\sum_{i,j=1}^{d}%
\frac{\partial }{\partial v_{i}}\left\{
\int_{\mathbb{R}^{d}}a_{ij}\left( v-v_{\ast }\right) \left[ f\left(t, v_{\ast }\right) \frac{%
\partial f}{\partial v_{j}}\left( t,v\right) -f\left( t,v\right) \frac{%
\partial f}{\partial v_{\ast j}}\left(t, v_{\ast }\right) \right]dv_{\ast
} \right\},
\end{equation}
with $a_{ij}(v):=(\sigma\sigma^*)_{ij}(v)=|v|^2\delta_{ij}-v_iv_j$
and $b_i(v)=\nabla \cdot a_{i\cdot}(v)$. The equations
\eqref{eq:Landaueq} model collisions of particles in a plasma and
can  be obtained as limit of the
 Boltzmann equations when collisions
become grazing, see Funaki \cite{Funaki:85},
 Goudon \cite{Goudon:97}, Villani \cite{Villani:98} \cite{Villani2:98} and Gu\'erin-M\'el\'eard \cite{Guerin:04}.

\smallskip

In this  work, we shall prove the convergence in law  of an easily
simulable mean field interacting particle system towards the
nonlinear process (\ref{eq:nonlinproc}) at an explicit pathwise
rate. This problem is of great interest in order to construct a
tractable simulation algorithm for the law $P_t$ and thus, in
particular, for solutions $f$ of equation \eqref{eq:Landaueq}. To
our knowledge, there is no result on convergence rates of the
deterministic numerical methods used at present for the Landau
equation, which are reviewed in \cite{Cordier:01}. The interest of
our approach is that it is based on the diffusive nature of the
equation, and that it addresses a large class of nonlinear
processes. The fact that we want to deal with simulable systems
will necessitate a coupling  between finite dimensional and
infinite dimensional stochastic processes. We shall introduce a
coupling argument based on new results on measurability of the
optimal mass transportation problem.

\smallskip

We  consider a particle system  which is naturally related to the
nonlinear process. Indeed, notice that the diffusion matrix
associated with (\ref{eq:nonlinproc}) is defined on $\RR^d $ by
 \be\label{eq:a}
 a(x,P_t):=\int_{\RR^d}\sigma(x-y)\sigma^*(x-y) P_t(dy)=[(\sigma\sigma^*)*
 P_t](x).
\ee Thus, if in order to approximate the
 white noise driven stochastic differential equation  \eqref{eq:nonlinproc}, we heuristically
   replace $P_t$ in \eqref{eq:a} by
an empirical measure of $n\in \NN^*$ particles in $\RR^d$, we are
led to consider the following system driven by $n^2$ independent
Brownian motions $(B^{ik})$:
\begin{equation}\label{eq:systparticules}
X^{i,n}_t=X_0^i+\frac{1}{\sqrt{n}}\int_0^t \sum_{k=1}^n
\sigma(X^{i,n}_s-X^{k,n}_s) dB_s^{ik}+\frac{1}{n}\int_0^t
\sum_{k=1}^n b(X^{i,n}_s-X^{k,n}_s) ds,~ i=1,\dots,n.
\end{equation}
To be more precise, if $\mu^n_t=\frac{1}{n}\sum_{i=1}^n
\delta_{X^{i,n}_t} $ is the empirical measure of the system, the
mappings
\begin{equation}\label{eq:martmeasyst}
f(t,\omega, x)\mapsto \frac{1}{\sqrt{n}}\int_0^t \sum_{k=1}^n
f(s,\omega,X^{k,n}_s) dB^{ik}_s, i=1,\dots,n,
\end{equation}
 define (for
suitably measurable functions $f$)  orthogonal martingale measures
in the sense of Walsh \cite{Walsh:84}, with covariance measure
$\mu^n_t\otimes dt$.

\smallskip
By adapting techniques  of M\'el\'eard-Roelly \cite{Meleard:88}
based on martingale problems, one can show propagation of chaos
for system \eqref{eq:systparticules} with as limit the process
\eqref{eq:nonlinproc}.  This says in particular that the
covariance measure of \eqref{eq:martmeasyst} converges in law to
$P_t\otimes dt$ when $n$ goes to infinity. But in turn, the
arguments of \cite{Meleard:88} do not give any information about
speed of convergence.

To estimate the distance between the law of the particles and the
law of the nonlinear process, we need to construct a significant
coupling between finitely many Brownian motions and the  white
noises processes. This problem is much more subtle than in the
McKean-Vlasov model (cf. Sznitman \cite{Szn} or M\'el\'eard
\cite{Mel}), where each particle is coupled with a limiting
process through a single Brownian motion that drives them both.
The well known $\frac{1}{\sqrt{n}}-$ convergence rate in that
model is consequence of the standard  $L^2$-law of large numbers
in $\RR^d$
 and of the fact that the diffusion and drift coefficients of the nonlinear process depend
linearly on the limiting law through  expectations with respect to
it. In the present Landau model, we have to deal with the
space-time random fields \eqref{eq:martmeasyst}, which have
fluctuations of constant order in $n$. This is also reflected  in
the fact that it is the  {\it squared} diffusion matrix of
\eqref{eq:nonlinproc}, that depends linearly  on $P_t$ (see
\eqref{eq:a}). It is hence not clear where a convergence rate can
be deduced from.

\medskip

Let $X^i$,$i=1,\dots,n$ be $n$ independent copies of the nonlinear
process in some probability space, and $\nu^n_t$  their empirical
measure at time $t$ (observe that it samples $P_t$). We shall
construct particles
  (\ref{eq:systparticules})  on the same probability space, in
  such way that they will converge  pathwise in
  $L^2$ on finite time intervals,
  at the same rate at which the Wasserstein distance $W_2$ between  $P_t$ and $\nu^n_t$ goes to $0$.
   Let us state our main
result on the process \eqref{eq:nonlinproc}:

\begin{theo}\label{main}
Let  $n\in \NN$ and assume usual  Lipschitz hypothesis on $\sigma$
and $b$, and that the law $P_0$ of $X_0^i$  has finite second
order moment.
 Assume moreover that $P_t$ has a density with respect to Lebesgue measure for each
 $t>0$.

Then, in the same probability space as $(X^1,\dots,X^n)$ there
exist independent standard Brownian motions $(B^{ik})_{1 \leq i,k
\leq n}$
 such that the particle
system $(X^{i,n})_{i=1}^n$ defined  in (\ref{eq:systparticules})
satisfies
\begin{equation*}
E\left(\sup_{t\in [0,T]}|X^{i,n}_t-X^i_t|^2\right)\leq C\exp(C'T)
\int_0^T E (W_2^2(\nu_s^n,P_s))ds
\end{equation*}
for constants $C,C'$ that do not depend on $n$.
\end{theo}

Thanks to available convergence  results for empirical measures of
i.i.d samples (see e.g. \cite{Rachev:98}), Theorem \ref{main} will
allow us to obtain, under some additional moment assumptions on
$P_0$, the speed of convergence $n^{{-2\over d+4}}$ for the
pathwise law of the system (see Corollary \ref{coro:weakrate}). We
remark that the absolute continuity condition of Theorem
\ref{main} can be obtained under non-degeneracy of the matrix
$\sigma\sigma^*$ by using for instance Malliavin calculus
\cite{Nualart:95}; it is also true for the specific coefficients
of the Landau equation \eqref{eq:Landaueq} despite their
degeneracy, and for some generalizations  (see Gu\'erin
\cite{Guerin:02}).

\smallskip

The proof of Theorem \ref{main} relies on new  results on the
optimal  mass transportation problem. For general background on
the theory of
 mass transportation, we refer to Villani
\cite{Villani:03}. Recall that if $\mu$ and $\nu$ are probability
measures in $\RR^d$ with finite second moment, the first of them
having a density, then the optimal mass transportation problem
with quadratic cost between $\mu$ and $\nu$ has a  unique
solution, which is a probability measure on $\RR^{2d}$ of the form
$\pi(dx,dy)= \mu(dx)\delta_{T(x)}(dy)$ . The so-called Brenier or
optimal transport map $T(x)$ is ($\mu$ a.s. equal to) the gradient
of some convex function in $\RR^d$, and pushes forward $\mu$ to
$\nu$.

\smallskip

Let now $W_P^i$  be the white noise process driving the $i$-th
nonlinear process $X^i$. The key idea in Theorem \ref{main} will
be to construct  Brownian motions
 $(B^{ik})_{k=1\dots n}$    in an  ``optimal'' pathwise way from  $W_P^i$.
 Heuristically, this will
  consist in pushing forward
the martingale measure $W_P^i$ through the Brenier maps
$T^{t,\omega,n}(x)$ realizing the optimal transport between $P_t$
and $\nu^n_t(\omega)$ (this is the reason for the absolute
continuity assumption on $P_t$). But to give such a construction a
rigorous sense, we must make sure that we can compute stochastic
integrals of $T^{t,\omega,n}(x)$ with respect to $W_P^i(dx,dt)$.
From the basic definition of stochastic integration with respect
to space-time white noise (cf. \cite{Walsh:84}),  this requires
the existence of a measurable version of $(t,\omega,x)\mapsto
T^{t,\omega,n}(x)$ being moreover predictable in $(t,\omega)$. A
striking problem then is that no available result in the mass
transportation theory can provide any information about joint
measurability properties of the optimal transport map, with
respect to the space variable and some parameter making the
marginals vary. Nevertheless, we will show that a suitable
``predictable transportation process'' exists:
\begin{theo}\label{theo:transportprocess}
There exists a measurable process $(t,\omega,x)\mapsto
T^{n}(t,\omega,x)$ that is predictable  in $(t,\omega)$ with
respect to the filtration associated to $(W^1_P,\dots,W^n_P)$ and
$(X_0^1,\dots,X^n_0)$, and such that for $dt \otimes \PP(d\omega)
\mbox{ almost every }(t,\omega),$
$$T^{n}(t,\omega,x)=T^{t,\omega,n}(x)\quad P_t(dx)\mbox{-almost surely.} $$
\end{theo}

  This
statement is consequence of a  general abstract result about
``measurability'' of the mass transportation problem. To be more
explicit, recall that the optimality of a transfert plan $\pi$ is
determined by its support (it is equivalent to the support being
cyclically monotone, see McCann \cite{McCann:95} or Villani
\cite{Villani:03}). On the other hand,  without assumptions
(besides moments) on the marginals $\mu$ and $\nu$, the solution
$\pi$ of the mass transportation problem may not be unique.  A
basic question then is how to formulate, in a general setting, the
adequate property of ``measurability'' of the solution(s) $\pi$
with respect to the data $(\mu,\nu)$. As we shall see, the natural
formulation  requires to introduce notions and techniques from
set-valued analysis. Then, we shall  prove the following

\begin{theo}\label{theo:measurable}
Let ${\cal P}_2(\RR^d)$ be the space of Borel probability measures
in $\RR^2$ with finite second order moment, endowed with the
Wasserstein distance and its Borel $\sigma-$field. Denote by
$\Pi^*(\mu,\nu)$ the set of solutions of the   mass transportation
problem with quadratic cost associated with $(\mu,\nu)\in ({\cal
P}_2(\RR^d))^2$. The function assigning to $(\mu,\nu)$ the  {\bf
set} of $\RR^{2d}$:
$$\bigcup_{ \pi \in \Pi^*(\mu,\nu)} supp (\pi),$$  is measurable in the sense of
set-valued mappings.
\end{theo}
In particular, this  ensures that if $\mu_{\lambda}$ and
$\nu_{\lambda}$ vary in a measurable way with respect to some
parameter $\lambda$, so that  in each of the associated optimal
transportation problems uniqueness holds, then the support of the
solution $\pi_{\lambda}$ also
 ``varies'' in a measurable way. This will be the key to our
 results.

\medskip

 The rest of this work is organized as follows. In Section
2 we review the Wasserstein distance and the mass transportation
problem with quadratic cost in $\RR^d$ (in particular  the
characterization of its minimizers). In Section 3 we prove Theorem
\ref{theo:measurable} and a consequence needed to prove Theorem
\ref{main}. In Section 4,   we state some properties about process
\eqref{eq:nonlinproc} and we heuristically describe our coupling
between space-time white noises and Brownian motions. In Section 5
we construct the ``predictable transportation process'' of Theorem
\ref{theo:transportprocess} needed to
  rigorously define the coupling. Section 6 is devoted to
 complete the proof of
 Theorem \ref{main} and to obtain explicit convergence rates.

\section{The mass transportation problem  with quadratic cost in $\RR^d$ and the Wasserstein
distance}

We denote the space of Borel probability measures in $\RR^d$ by
${\cal P}(\RR^d)$, and by ${\cal P}_2(\RR^d)$ the subspace of
probability measures having finite second order moment.

Given $\pi \in {\cal P}_2(\RR^{2d})$, we respectively denote by
$\pi_1$ and $\pi_2$ its first and second marginals on $\RR^d$. On
the other hand, for any two probability measures $\mu,\nu\in {\cal
P}_2(\RR^d)$ and $\pi \in {\cal P}_2(\RR^{2d})$, we write
$$\pi<^{\mu}_{\nu}$$
if $\pi_1=\mu$ and $\pi_2=\nu$.
Such $\pi$ is refereed to as a ``transfert plan'' between $\mu$
and $\nu$.

\begin{defi}
The Wassertein distance $W_2$ on ${\cal P}_2(\RR^d)$ is defined by
$$W_2^2(\mu,\nu):=\inf_{\pi<^{\mu}_{\nu}} \int_{\RR^2} |x-y|^2 \pi(dx,dy).$$
\end{defi}

Then, $({\cal P}_2(\RR^d),W_2)$ is a Polish space, see e.g. Rachev
and R\"uschendorf \cite{Rachev:98}. The topology is stronger that
the usual weak topology. More precisely, one has the following
result (see for instance Villani, \cite{Villani:03} Theorem 7.12)

\begin{theo}\label{theo:Wasser} Let $\mu^n,\mu\in {\cal
P}(\RR^d)$. The following are then equivalent:
\begin{itemize}
\item[i)] $W_2(\mu^n,\mu)\to 0$ when $n\to \infty$.
\item[ii)] $\mu^n$ converges weakly to $\mu$ and
$$ \int_{\RR^d} |x|^2 \mu^n(dx)\to \int_{\RR^d} |x|^2 \mu(dx).$$
\item[iii)] We have
$$ \int_{\RR^d} \varphi(x) \mu^n(dx)\to \int_{\RR^d} \varphi(x) \mu(dx)$$
for all continuous function $\varphi:\RR^d\to \RR$ such that
$|\varphi(x)|\leq C(1+|x|^2)$ for some $C\in \RR$.
\end{itemize}
\end{theo}

We shall denote by $L$ the mapping $L:{\cal P}_2(\RR^{2d})\to\RR$
defined by
$$L(\pi)=\int_{\RR^2} |x-y|^2 \pi(dx,dy).$$

\begin{rem}\label{rem:contL}
It is not hard to check that $L$ is lower semi continuous (l.s.c)
for the weak topology.  Moreover,  $L$ is continuous for the
Wasserstein topology in ${\cal P}_2(\RR^{2d})$ by part {\it iii)}
of Theorem \ref{theo:Wasser}.
\end{rem}

\medskip

Fix now $\mu,\nu \in {\cal P}_2(\RR^d)$, and denote by
$\Pi^*(\mu,\nu)$ the subset of ${\cal P}_2(\RR^{2d})$ of
minimizers of the {\it Monge-Kantorovich transportation problem}
with quadratic cost for the pair of marginals $(\mu,\nu)$ . That
is,
$$\Pi^*(\mu,\nu):=argmin_{ \pi<^{\mu}_{\nu}} L(\pi).$$

It is well known that $\Pi^*(\mu,\nu)$ is non-empty. Indeed, it is
not hard to see that for the weak topology, $\{\pi \in {\cal
P}_2(\RR^{2d}): \pi<^{\mu}_{\nu}\}$ is  a compact set,  and the
lower semi-continuity of $L$ implies the existence of minimizers
(see e.g.  \cite{Villani:03} Chapter 1 for details).

\bigskip

We shall next recall the characterization of minimizers of the
transportation problem with quadratic cost. We need the notion of
sub-differential of a convex function:
\begin{defi}
Let $\varphi:A\subset \RR^d \to ]-\infty,\infty]$ be  a proper
(i.e. $\varphi\not \equiv +\infty$) lower semi-continuous (l.s.c)
convex function. The sub-differential  of $\varphi$ at $x$ is
$$ \partial \varphi(x)=\{y\in \RR^d: \varphi(z)\geq \varphi(x) +
\langle y,z-x\rangle , \forall z \in \RR^d\}.$$

Elements of $\partial \varphi(x)$ are called sub-gradients of
$\varphi$ at point $x$. The graph of $\partial \varphi$ is
$$Gr(\partial \varphi)=\{(x,y)\in \RR^{2d}: y\in \partial
\varphi(x)\}$$ and it is a closed set.

\end{defi}
Recall that $\varphi$ is differentiable at $x$ if and only if
$\partial \varphi(x)$ is a singleton (in which case $\partial
\varphi(x)=\{\nabla \varphi(x)\}$). Also, the set   $\{x\in \RR^d:
\varphi\mbox{ is differentiable at }x\}$ is borelian, see e.g. McCann
\cite{McCann:95}.

\smallskip

We next summarize  results in pioneer works in this domain,
Knott-Smith \cite{Knott:84}, Brenier \cite{Brenier:91} and McCann
\cite{McCann:95}, Rachev and R\"uschendorf \cite{Rachev:98}. See
also Villani \cite{Villani:03} for a complete discussion on these
questions, proofs and background.

\begin{theo}\label{theo:charactmin}
Let $\mu,\nu \in {\cal P}(\RR^d)$ and $\pi<^{\mu}_{\nu}$ be a
transfert plan. We have
\begin{itemize}
\item[a)] $\pi \in \Pi^*(\mu,\nu)$ if and only if there exists a
proper l.s.c. convex function $\varphi$ such that $$
supp(\pi)\subset Gr(\partial \varphi)$$ or, equivalently
$$\pi(\{(x,y)\in \RR^2: y\in
\partial \varphi(x)\})=1.$$
\item[b)] Assume that   $\mu$  does not charge sets of Hausdorff
dimension less or equal than $d-1$ and that $\pi \in
\Pi^*(\mu,\nu)$. Then,
\begin{itemize}
\item[i)] the set $\{x\in \RR^d: \varphi\mbox{ is not differentiable
at }x\}$ has null $\mu$-measure.
\item[ii)] We have
$$\pi(dx,dy)=\mu(dx)\otimes \delta_{\nabla \varphi(x)}(dy).$$
\item[ii)] If $T$ is a measurable mapping such that $\pi(dx,dy)=\mu(dx)\otimes
\delta_{T(x)}(dy)$, then $T(x)=\nabla \varphi(x)$ ,
$\mu(dx)-a.s.$.
\item[iii)]  $\pi \in \Pi^*(\mu,\nu)$ is unique.
\end{itemize}
\end{itemize}
\end{theo}

This result will be  useful later  in the particular case when
the measure $\mu$ is absolutely continuous with respect to
Lebesgue measure.

\section{Measurability of the mass transportation problem}

 We now introduce the basic notions on ``multi-applications'' or
``set-valued mappings'' that we  need to prove Theorem
\ref{theo:measurable}. For general background, we refer the reader
to  Appendix A in Rockafellar and Wets \cite{Rocka:98}.

\begin{defi}
Let $X,Y$ be two sets.
\begin{itemize}
\item[i)] A function $S$ on $X$ taking values in
the set of subsets of $Y$ is called a set-valued mapping or
multi-application. We write $S:X \rightrightarrows Y$.

\item[ii)] For any $A\subset Y$, the inverse image of $A$ through $S$ is
the set
$$S^{-1}(A):=\{x\in X: S(x)\cap A \not = \emptyset\}.$$

\item[iii)] If $(X,{\cal A})$ is a measurable space and $(Y,\Theta)$
a topological space, we say that $S:X \rightrightarrows Y$ is
measurable if for all $\theta \in \Theta$,
$$S^{-1}(\theta)\in {\cal A}.$$
(Of course, if $S(x)=\{s(x)\}$ is singleton for all $x$,
measurability of $S$ is equivalent to that of $s$. )
\end{itemize}

\end{defi}

Consider  ${\cal P}_2(\RR^d)$  endowed with the Wasserstein
distance and the Borel $\sigma-$field. We  define a  set-valued
mapping $$\Psi:({\cal P}_2(\RR^d))^2 \rightrightarrows\RR^{2d}$$
by
$$ \Psi(\mu,\nu):=\{(x,y): \exists \pi \in \Pi^*(\mu,\nu) ~s.t.~
(x,y)\in supp (\pi)\}.$$

\medskip

 Our goal is to prove that $ \Psi$ is measurable.  We shall
need some further notions on set-valued mappings.

\begin{defi}
Let $X$ be a set, and $(Y,\Xi)$ and $(Z,\Theta)$ be topological
spaces.

\begin{itemize}
\item[i)] A  set-valued mapping $S:X\rightrightarrows Y$ is
closed-valued if for all $x\in X$, $S(x)$ is a closed set of
$(Y,\Xi)$.

\item[ii)] A set-valued mapping $ U:Y\rightrightarrows Z$ is inner
semicontinous (i.s.c) if for all $\theta \in \Theta$,
$$S^{-1}(\theta)\in \Xi$$

\end{itemize}

\end{defi}

The following results can be found in Appendix A of
\cite{Rocka:98}, in the case of set-valued mappings in $\RR^d$.
For completeness we provide proofs in a more general context.

\begin{lem}\label{lem:propsmultiapplic}
 Let $(X,{\cal A})$ be a measurable space and $(Y,\Xi)$ a
 topological  space.
\begin{itemize}
\item[i)] $S:X \rightrightarrows Y$ is measurable if and only
if the closed-valued mapping  $x\rightrightarrows Cl(S(x))$ is
measurable, where $Cl(S(x))$ is the topological closure of the set
$S(x)$.
\item[ii)] Assume that $(Y,d)$
 is a separable metric space and that $S:X\rightrightarrows Y$ is
closed-valued. Then,  $S$ is measurable if and only if for all
closed set $F$ of ~ $Y$,
$$S^{-1}(F)\in {\cal A}.$$

\item[iii)] Let $(Y,\Xi)$ and $(Z,\Theta)$ be topological spaces,
$S:X\rightrightarrows Y$ be measurable and $U:Y\rightrightarrows
Z$ be i.s.c. Then, the multi-application $U\circ S:X
\rightrightarrows Z$, defined by
$$U\circ S(x):=\bigcup_{y\in S(x)} U(y)$$
is measurable.

\end{itemize}
\end{lem}

{\bf Proof} {\it i)} For any open set $\theta\in \Xi$, $ S(x)\cap
\theta \not = \emptyset\mbox{ if and only if }Cl(S(x))\cap \theta
\not = \emptyset.$

{\it ii)}  ``Only if'' part: since $Y$ is a metric space, we use
that every closed set $F$ is the intersection of some countable
collection of open sets $(\theta_n)$. Therefore,
$$\{x\in X:S(x)\cap F\not = \emptyset \}=\bigcap_{n\in \NN} \{x\in X:S(x)\cap \theta_n\not = \emptyset
\}\in {\cal A}.$$

 ``If'' part:  $(Y,d)$ being separable, we can express every open
set $\theta$ as the union of some countable collection $(B_n)$ of
closed balls. We then have that
$$\{x\in X:S(x)\cap \theta\not = \emptyset \}=\bigcup_{n\in \NN} \{x\in X:S(x)\cap B_n\not = \emptyset
\}\in {\cal A}.$$

{\it iii)} Straightforward:
\begin{eqnarray*}
(U\circ S)^{-1}(\theta)&= & \{x\in X: \left(\cup_{ y\in S(x)}
U(y)\right)\cap \theta \not = \emptyset\}=  \{x\in X: ~\exists
y\in S(x) ~s.t. ~ U(y)\cap \theta \not = \emptyset\}\\
  &= & \{x\in X: S(x)
\cap (U^{-1}( \theta)) \not = \emptyset\}.
\end{eqnarray*}

The function $U$ being i.s.c.,  $U^{-1}( \theta)$ belongs to
$\Xi$, which allows us to conclude.

\Box

\medskip

Now we can proceed to the

\smallskip

{\bf Proof of Theorem \ref{theo:measurable}}

 We observe first that
$\Psi(\mu,\nu)=U\circ S(\mu,\nu)$, where  $S$ and $U$ are the set
valued mappings respectively defined by
$$(\mu,\nu)\rightrightarrows S(\mu,\nu):= \Pi^*(\mu,\nu)$$ and
 $U:{\cal P}_2(\RR^{2d})\rightrightarrows \RR^d$ by
$$U(\pi):=supp(\pi)$$

   We will therefore  split the
proof in several parts:

{\it a) $S$ is a closed valued mapping }

First notice that  $\pi \mapsto \pi_i$ is continuous for the
Wasserstein topology. Indeed, $W_2(\pi^n,\pi)\to 0$ implies that
$\pi^n$ converges weakly to $\pi$, and then $\pi^n_i$ converges
weakly to $\pi_i$ for $i=1,2$. Moreover, we have $\int_{\RR^d}
|x|^2\pi_1^n(dx)=\int_{\RR^{2d}} |x|^2 \pi^n(dx,dy)\to
\int_{\RR^{2d}} |x|^2 \pi(dx,dy)=\int_{\RR^d} |x|^2\pi_1(dx)$ by
Theorem \ref{theo:Wasser}, and then the asserted continuity
follows.

Consequently, $ \pi \mapsto W_2(\pi_1,\pi_2) $ too is continuous.
Therefore,
$$\Pi^*(\mu,\nu) =\{\pi: \pi<^{\mu}_{\nu}\} \cap \{\pi: L(\pi)-W_2(\pi_1,\pi_2)
=0\}$$ is the intersection of two closed sets  ${\cal
P}_2(\RR^{2d})$.

\medskip

 {\it b) Inverse images through $S$ of closed sets  are
closed sets  }

Let $F \subset {\cal P}_2(\RR^d)$ be a closed set and
$(\mu^n,\nu^n)\in S^{-1}(F)$,\ $n\in \NN$ , be a sequence
converging to $(\mu,\nu)$ in $({\cal P}_2(\RR^d))^2$. Then,
$\mu^n\to \mu$ and $\nu^n\to \nu$ weakly, and $(\mu^n)$ and
$(\nu^n)$ are tight.

But since $(\mu^n,\nu^n)\in S^{-1}(F)$ for each $n$, there exists
$\pi_n$ s.t. $\pi^n<^{\mu^n}_{\nu^n}$,  and then $(\pi_n)$ too is
tight (by considering products of compact sets).

Let $(\pi^{n_k})$ be a weakly convergent subsequence with limit
$\pi$. Then, clearly $\pi<^{\mu}_{\nu}$. We will  prove that
$L(\pi)=W_2(\mu,\nu)$ and that $\pi\in F$, which will mean that
$(\mu,\nu)\in S^{-1}(F)$ and finish the proof.

We have
\begin{multline*}
\int_{\RR^{2d}} \left(|x|^2+|y|^2\right)
\pi^{n_k}(dx,dy)=\int_{\RR^{d}} |x|^2 \mu^{n_k}(dx)+\int_{\RR^{d}}
|y|^2 \nu^{n_k}(dy) \to \\
\int_{\RR^{d}} |x|^2 \mu(dx)+\int_{\RR^{d}} |y|^2 \nu(dy)=
\int_{\RR^{2d}} \left(|x|^2+|y|^2\right) \pi(dx,dy),
\end{multline*}
which implies that $W_2(\pi^n,\pi)\to 0$ and $\pi \in F$. Finally,
by the continuity of $\pi\mapsto L(\pi)- W_2(\pi_1,\pi_2)$ we get
that
$$0=L(\pi^{n_k})-
W_2(\pi^{n_k}_1,\pi^{n_k}_2)=L(\pi)-W_2(\mu,\nu).$$

\medskip

 {\it c) The mapping $U$ is  i.s.c. }

Let $\theta$ be an open set of $\RR^{2d}$. We must check that
\begin{equation*}
 \{\pi \in {\cal P}_2(\RR^{2d}): supp(\pi) \cap \theta \not =
\emptyset\} = \{\pi \in {\cal P}_2(\RR^{2d}): \pi(\theta)>0 \}
\end{equation*}
 is open, or equivalently, that
$$\{\pi \in {\cal P}_2(\RR^{2d}): \pi(\theta)=0 \}$$
is closed in $ {\cal P}_2(\RR^{2d})$. Assume that  $\pi,\pi^n\in
{\cal P}_2(\RR^{2d})$, with $\pi^n$ such that $\pi^n(\theta)=0 $
for all $n\in \NN$, and moreover that $W_2(\pi^n,\pi)\to 0$. Then
$\pi^n$ converges weakly to $\pi$, and so by the Portemanteau
theorem, we have
$$0=\liminf_n \pi^n(\theta)\geq \pi(\theta).$$

{\it d) Conclusion}

 By parts {\it a)} and {\it b)} and  Lemma \ref{lem:propsmultiapplic} {\it ii)}
 we get that $S$ is measurable. By {\it c)} and Lemma  \ref{lem:propsmultiapplic} {\it
 iii)} $U\circ S$ is measurable and the proof is finished.

\Box

\medskip

 The following corollary will be useful in the specific
setting needed to prove Theorem \ref{main}:

\begin{coro}\label{coro:useful}

Let $(E,\Sigma)$ be a measurable space,  and $ \lambda \in E
\mapsto (\mu_{\lambda},\nu_{\lambda})\in ({\cal P}_2(\RR^{d}))^2$
and
 $\xi: E\to \RR^d$ be measurable functions.
Then, the set
 $$\left\{(\lambda,x): (x,\xi(\lambda))\in
 Cl(\Psi)(\mu_{\lambda},\nu_{\lambda})\right\}$$ belongs to
$ \Sigma\otimes {\cal B}(\RR^d)$
\end{coro}

{\bf Proof} By Lemma \ref{lem:propsmultiapplic} {\it
 i)} and Theorem \ref{theo:measurable} we get that $Cl(\Psi)$ is measurable.
 Moreover, it is not hard to check that
the mapping
$$(\lambda,x)\rightrightarrows
Cl(\Psi)(\mu_{\lambda},\nu_{\lambda})-(x,\xi(\lambda))$$ is
measurable and closed-valued. Then, we just have to notice that

$$(x,\xi(\lambda))\in
Cl( \Psi)(\mu_{\lambda},\nu_{\lambda}) \mbox{ if and only if }
\left[
Cl(\Psi)(\mu_{\lambda},\nu_{\lambda})-(x,\xi(\lambda))\right]\cap
C\not = \emptyset$$ for the closed set $C=\{0\}$.

\Box

\section{A coupling between space-time white noise and Brownian motions  via optimal transport}

In all the sequel, we refer the reader to
 Walsh \cite{Walsh:84} for background on space-time white noise processes and  stochastic
integration with respect to martingale measures.

Assume  that  $\sigma:\RR^d\to \RR^{d\otimes d}$ and $b:\RR^d\to
\RR^d$ are Lipschitz continuous and with linear growth. Then, by
results of \cite{Funaki:84} or \cite{Guerin:03} we can construct
in some probability  space $(\Omega,{\cal F},\PP)$  a sequence
 $(X^i)_{i\in \NN}$  of independent copies of the nonlinear
processes,
\begin{equation}\label{eq:indprocnonlin}
X^i_t=X^i_0+\int_0^t \int_{\RR^d} \sigma(X^i_s-y)
W^i_P(dy,ds)+\int_0^t \int_{\RR^d} b(X^i_s-y) P_s(dy)ds,
\end{equation}
where the $W^i_P$ are independent space-time $\RR^d$-valued white
noises defined on $[0,\infty)\times \RR^d$. Each of the $d$
(independent) coordinates of $W^i_P$ has covariance measure
$P_t(dy)\otimes dt$, where $P_t$ is the law of $X_t$. The initial
conditions $(X_0^1,\dots,X_0^n,\dots)$ are independent and
identically distributed with law $P_0$, and independent of the
white noises. The pathwise law of $X^i$ is denoted by $P$, and it
is uniquely determined.

\smallskip

 Denote by ${\cal F}^n_t$ the
complete right continuous $\sigma$-field generated by
$$\{(W_P^1([0,s]\times A^1),\dots,W_P^n([0,s]\times A^n)): 0\leq
s\leq t, A^i\in {\cal B}(\RR^d)\}$$ and $(X_0^1,\dots,X_0^n).$
  We also denote by   $${\cal P}red^n$$
   the predictable field
generated by continuous $({\cal F}^n_t)$-adapted processes.

\medskip

In what follows, we fix a finite time horizon $T>0$. Under usual
Lipschitz assumptions on the coefficients, there is propagation of
the moments of the law $P_0$, as proved in Gu\'erin
\cite{Guerin:03}.
\begin{lem}\label{lem:momboud}
\label{moment}
 If $E(|X_0|^k)<\infty$ for some $k \geq 2$, then
\begin{equation*}
 E\left(\sup_{t\in [0,T]}|X_t|^k\right)<\infty.
\end{equation*}
The continuity of $X$ and the previous uniform bound  imply that
 $t\mapsto \int_{\RR^d} |x|^k
P_t(dx)$ is continuous.
\end{lem}

Throughout the sequel, the assumptions of Theorem \ref{main} on
$P_0$ and $P_t$ are enforced, in particular, the condition
$E(\sup_{t\in [0,T]}|X_t|^2)<\infty$ will hold by the previous
lemma.

\medskip We shall now present the main idea of the coupling we
introduce to prove Theorem \ref{main}. Basically, this consists in
constructing for each $n$,  $n^2$  Brownian motions in a pathwise
way, from the realizations of the $n$  white noises
$(W^1_P,\dots,W_P^n)$. The key for that will be to use the optimal
 transport maps between the marginal $P_t$ of the nonlinear
process and the empirical measures of samples of that law. More
precisely, write
\begin{equation*}
 \nu_t^n:=\frac{1}{n}\sum_{i=1}^n
\delta_{X^i_t}
\end{equation*}
and notice that for each  $\omega\in \Omega$, $(\nu_t^n, 0\leq
t\leq T)$ is an element of $C([0,T],{\cal P}_2(\RR^d))$. Thus, for
each $t\in [0,T]$, $n\in \NN$ and $\omega$, and we can consider
the optimal coupling problem with quadratic cost between
$\nu_t^n(\omega)$ and $P_t$,
\begin{equation*}
\inf_{\pi<^{P_t}_{\nu_t^n(\omega)}  } \left\{ \int_{\RR^d\times
\RR^d}  |x-y|^2 \pi(dx,dy)\right\}.
\end{equation*}

\medskip
By the assumption on $P_t$ and Theorem \ref{theo:charactmin}, the
following properties hold {\bf for each fixed pair} $(t,\omega)\in
]0,T]\times \Omega$:

\begin{lem}
\begin{itemize}
\item[a)] There exists a unique  $\pi^{t,\omega,n}$,
such that
\begin{equation*}
W_2^2(P_t,\nu_t^n(\omega))=\int_{\RR^d\times \RR^d}  |x-y|^2
\pi^{t,\omega,n}(dxdy).
\end{equation*}
\item[b)] There is a $P_t(dx)-a.e.$ unique measurable function
$T^{t,\omega,n}:\RR^d\to \RR^d$ such that
\begin{equation*}
\pi^{t,\omega,n} (dx,dy)=\delta_{T^{t,\omega,n}(x)}(dy)P_t(dx).
\end{equation*}
In particular, under $P_t(dx)$ the law  of $T^{t,\omega,n}(x)$ is
$\nu_t^n(\omega)$. \item[c)] We have
\begin{equation*}
W_2^2(P_t,\nu_t^n(\omega))=\int_{\RR^2} |x-T^{t,\omega,n}(x)|^2
P_t(dx).
\end{equation*}
\end{itemize}
\end{lem}
We would like to construct $n^2$ independent Brownian motions by
``transporting'' the $n$ independent white noises
$(W_P^1,\dots,W_P^n)$ through the transport mappings
$T^{s,\omega,n}(x)$. As pointed out in the introduction, to do so
we must at least be able to define stochastic integrals of
functions of the form $(t,\omega,x)\mapsto f(T^{t,\omega,n}(x))$,
with respect to the white noise processes.
 The existence of  a version $T^n(t,\omega,x)$ of $ T^{t,\omega,n}(x)$
 having good enough properties, will be established in next section,
 when we shall prove Theorem \ref{theo:transportprocess}.

\smallskip Before doing so, we observe that if  Theorem
\ref{theo:transportprocess} holds, then the following processes
$B_t^{ik}=B_t^{ik,n}$ will be well defined from
(\ref{eq:indprocnonlin}).

\begin{prop}\label{prop:Brms}
For each $n\in \NN^*$, define
\begin{equation}\label{eq:lesbrowns} B_t^{ik,n}(\omega):=\sqrt{n}\int_0^t \int_{\RR^d}
\1_{\{T^n(s,\omega,x) =X^{k}_s(\omega)\}}W_P^i(dx,ds),\quad
i,k=1\dots n
\end{equation}
 Then, $(B^{ik,n})_{1 \leq i,k \leq n}$ are $n^2$
independent standard Brownian motions in $\RR^d$.
\end{prop}

These are  the right Brownian motions we need to construct
\eqref{eq:systparticules}. The proof of Proposition
\ref{prop:Brms} will given in Section 6.

\section{Construction of the predictable ``transport process''}

Our goal now in this section is to show that for each $n \in
\mathbb{N}^*$, there exists a process $(t,\omega,x)\mapsto
T^{n}(t,\omega,x)$ defined $\PP(d\omega)\otimes dt \otimes
P_t(dx)$-almost everywhere,  which is measurable with respect to
${\cal P}red^n\otimes {\cal B}(\RR^d)$, and such that:
$$\mbox{ for }dt \otimes \PP(d\omega) \mbox{ almost every }(t,\omega),$$
$$T^{n}(t,\omega,x)=T^{t,\omega,n}(x)\quad P_t(dx)\mbox{-almost surely .} $$

 Since $(X^1,\dots,X^n)$ are independent copies of the
 nonlinear process and $P_t=law(X^i_t)$ has a density with respect to Lebesgue measure,  for each $t\in
 [0,T]$ we have that
 $$\PP(\exists i\not =j: X^i_t=X^j_t)=0.$$
Notice also that for fixed $(i,j)$ with  $i\not =j$  the following
set
$$\{(t,\omega):X^i_t(\omega)=X^j_t(\omega)\}$$ belongs to ${\cal
P}red^n$ since $(t,\omega)\mapsto |X^i_t(\omega)-X^j_t(\omega)|$
is adapted and  continuous in $t$.

 By Fubini's  theorem we then see that
$$ \int_{[0,T]\times \Omega} \1_{\{
 X^i_{t}=X^j_{t} \}}(t,\omega)\PP(d\omega)\otimes dt=0$$

\begin{rem}\label{rem:nullset}
Consequently,  there is a predictable set of $[0,T]\times
\Omega$,
$$\Omega'_T\in {\cal P}red^n$$ of full $\PP(d\omega)\otimes
dt$-measure and  such that
$$\mbox{ for all }(t,\omega)\in \Omega'_T,\quad
X^i_t(\omega)\not =X^j_t(\omega)\mbox{ for all }i,j\in
\{1,\dots,n\}.$$
\end{rem}
Let us denote by $({\cal P}red^n)'$   the $\sigma-$field ${\cal
P}red^n$ restricted to $\Omega'_T$.

Recall that for each $(t,\omega)$,  the set of solutions
$\Pi^*\left(P_t,\frac{1}{n}\sum_{i=1}^n
\delta_{X^i_t}(\omega)\right)$
 of the optimal transport problem between $P_t$ and $\frac{1}{n}\sum_{i=1}^n
 \delta_{X^i_t}(\omega)$ is a singleton that we have denoted by
$\pi^{t,\omega,n}$.

Let us define now the sets

$$A^{i,n}:=\left\{(t,\omega,x)\in \Omega'_T\times \RR^d:
(x,X^i_t(\omega))\in supp(\pi^{t,\omega,n})\right\},\quad
i=1,\dots,n.$$

The sets $A^{i,n}$ are predictable, as proved in the following
lemma.

\begin{lem} We have $A^{i,n}\in ({\cal P}red^n)'\otimes
{\cal B}(\RR^d)$.
\end{lem}

{\bf Proof} Observe that  the deterministic process
$(t,\omega)\mapsto P_t\in {\cal P}_2(\RR^d)$ is ${\cal
P}red^n$-measurable. Indeed, if  $(f_n)_{n\in \NN\backslash
\{0,1\}}$ is a countable dense subset of  the space of continuous
functions in $\RR^d$ with compact support, and
$f_0(x)=1$,$f_1(x)=|x|^2$, then the topology of ${\cal
P}_2(\RR^d)$ is generated by the real mappings $m\mapsto \int
f_n(x) m(dx)$.  It is therefore enough that $(t,\omega)\to \int
f_n(x) P_t(dx)$ be ${\cal P}red^n$-measurable, which is clear
since $t\mapsto P_t$ is continuous.

 Next we will apply   Corollary
\ref{coro:useful} to the measurable space
$$(E,\Sigma)=( \Omega'_T, ({\cal P}red^n)'),$$
$\lambda=(t,\omega)$, and the $({\cal P}red^n)'$-measurable
functions given by
$$(t,\omega)\to \left(P_t,\frac{1}{n}\sum_{j=1}^n \delta_{X^j_t}(\omega)\right)\in {\cal P}_2(\RR^{2d})~\mbox{
and }~ (t,\omega)\mapsto \xi^i(t,\omega)=X^i_t(\omega)\in \RR^d.$$
 For each $(t,\omega)\in
\Omega'_T$, with $\Psi$ denoting the multi-application defined in
Theorem \ref{theo:measurable}, we simply have in the current
setting that
$$\Psi\left(P_t,\frac{1}{n}\sum_{j=1}^n
\delta_{X^j_t}(\omega)\right)=supp(\pi^{t,\omega,n}).$$ Corollary
\ref{coro:useful} implies the result.

\Box

\medskip

Recall from basic measure theory that if $E_1$ and $E_2$ are
measurable spaces and $A\subseteq E_1\times E_2 $ is an element of
their product $\sigma-$field, then, for each $\lambda_1\in E_1$,
the {\it fiber} of $A$ at $\lambda_1$ is the set
$$[A]_{\lambda_1}:=\{\lambda \in E_2 : (\lambda_1,\lambda)\in A\},$$
and it is always measurable in $E_2$.

We can now proceed to the

\smallskip

 {\bf Proof of Theorem \ref{theo:transportprocess}:}
\medskip

We split the proof in several parts.

\medskip

 {\it a) The sets
$A^{i,n},i=1\dots n$ form a partition of $\Omega'_T\times \RR^d$
up to $\PP(d\omega)\otimes dt \otimes P_t(dx)$-null sets.}

For  $i\not = j$ write
\begin{eqnarray*}
A^{ij,n}&:=&\{(t,\omega,x)\in (\Omega_T'\times \RR^d):
(x,X^i_t(\omega))\in supp(\pi^{t,\omega,n})\mbox{ and }
(x,X^j_t(\omega))\in
supp(\pi^{t,\omega,n})\}\\
&=& A^{i,n}\cap A^{j,n}, \end{eqnarray*}
 and denote by
$[A^{ij,n}]_{(t,\omega)}:=\{x\in \RR^d: (t,\omega,x) \in
A^{ij,n}\}\in {\cal B}(\RR^d)$ the fiber of $A^{ij,n}$ at
$(t,\omega)\in \Omega'_T$. Then, we have
\begin{equation*}
\begin{split}
P_t([A^{ij,n}]_{(t,\omega)})= & ~ P_t(\{x\in
\RR^{d}:(x,X^i_t(\omega)),(x,X^j_t(\omega))\in
supp(\pi^{t,\omega,n})\}) \\
\leq & ~ P_t(\{x\in \RR^{d}:X^i_t(\omega),X^j_t(\omega)\in
\partial
\varphi^{t,\omega,n}(x)\}) ,\\
  \end{split}
\end{equation*}
where  $\varphi^{t,\omega,n}$ is a proper l.s.c. convex function
given by Theorem \ref{theo:charactmin} {\it a)}. But since
$(t,\omega)\in \Omega'_T$, we have $X^i_t(\omega)\not
=X^j_t(\omega)$, and so
$$X^i_t(\omega),X^j_t(\omega)\in \partial
\varphi^{t,\omega,n}(x) \Longrightarrow \varphi \mbox{ is not
differentiable in }x.$$ We obtain by Theorem \ref{theo:charactmin}
{\it b)} that $P_t([A^{ij,n}]_{(t,\omega)})=0,$ and then
$$E\left(\int_{[0,T]\times \Omega\times \RR^d}{\bf 1}_{ A^{i,n}\cap
A^{j,n}}(t,\omega,x)P_t(dx)dt\right)=0.
$$

\medskip

On the other hand,  since $T^{t,\omega,n}(x)\in
\{X_t^1(\omega),\dots,X_t^n(\omega)\}$  $P_t(dx)$ a.s., we have
for all $(t,\omega)$ that
\begin{equation*}
\begin{split}
P_t\left(\left[\left(\bigcup_{i=1}^n
A^{i,n}\right)^c\right]_{(t,\omega)}\right) & =P_t(\{x\in \RR^d:
\mbox{ for all }i=1,\dots,n,
 (x,X^i_t(\omega))\not\in supp(\pi^{t,\omega,n})\}) \\
 & \leq  P_t(\{x\in \RR^d:
 (x,T^{t,\omega,n}(x))\not\in supp(\pi^{t,\omega,n})\}) \\
& = \pi^{t,\omega,n}( supp(\pi^{t,\omega,n})^c )\\
& = 0 \\
\end{split}
\end{equation*}

Defining the set
$$\tilde{\Omega}_T:=(\Omega'_T\times \RR^d)\bigcap \left( \bigcup_{i=1}^n
A^{i,n} \backslash \left(\bigcup_{k\not = j}
A^{kj,n}\right)\right) \in {\cal P}red^n\otimes {\cal B}(\RR^d)$$
we deduce  that
$P_t\left([\tilde{\Omega}_T^c]_{(t,\omega)}\right)=0$ for all
$(t,\omega)\in \Omega'_T$. Therefore,
$$\EE\left( \int_0^T \int_{\RR^d} \1_{\tilde{\Omega}_T^c}(t,\omega,x)
P_t(dx) dt\right)=\int_{\Omega'_T}
P_t\left([\tilde{\Omega}_T^c]_{(t,\omega)}\right) ~ dt\otimes
\PP(d\omega)=0,$$ and so $\tilde{\Omega}_T$ has full
  $\PP(d\omega)\otimes dt
\otimes P_t(dx)-$measure. This proves assertion {\it a)}.

\medskip

   We
can now define a ${\cal P}red^n\otimes {\cal B}(\RR^d)-$measurable
function by
\begin{equation}\label{predtrans}
T^n(t,\omega,x):=\sum_{i=1}^n \1_{A^{i,n}\cap
\tilde{\Omega}_T}(t,\omega,x) X^i_t(\omega).
\end{equation}

\medskip

{\it b) For  $\PP(d\omega)\otimes dt$ almost every $(t,\omega)$,
$T^n(t,\omega,x)=T^{t,\omega,n}(x)$ holds $P_t(dx)$ almost
surely.}

By Theorem \ref{theo:charactmin}, {\it b)}, this is equivalent to
prove that
\begin{equation*}\label{eq:T=T}
\pi^{t,\omega,n}(dx,dy)=P_t(dx)\otimes
\delta_{T^n(t,\omega,x)}(dy)\quad \PP(d\omega)\otimes dt-a.e.
\end{equation*}
We fix now $(t,\omega)\in \Omega'_T$ and $C,D\in {\cal B}(\RR^d)$.

We have by definition of  $T^{t,\omega,n}$ that
\begin{equation*}
\begin{split}
\pi^{t,\omega,n}(C\times D)= & \int_{\RR^d}
\1_C(x)\1_D(T^{t,\omega,n}(x))P_t(dx) \\
= & \int_{\RR^d} \1_{C\cap
[\tilde{\Omega}_T]_{(t,\omega)}}(x)\1_D(T^{t,\omega,n}(x))P_t(dx),\\
\end{split}
\end{equation*}
the latter because
$P_t\left([\tilde{\Omega}^c]_{(t,\omega)}\right)=0$.
 Notice that
on the other hand, by definition of $A^{i,n}, \tilde{\Omega}_T$
and $T^n$, for all $(t,\omega,x)\in A^{i,n}\cap \tilde{\Omega}_T$
we have that
$$\{y:(x,y)\in supp(\pi^{t,\omega,n})\}=\{X^i_t(\omega)\}=\{T^n(t,\omega,x)\}.$$
This implies that $\tilde{\Omega}_T\subset \{(t,\omega,x)\in
\Omega'_T\times\RR^d: \{y:(x,y)\in supp(\pi^{t,\omega,n})\} \mbox{
is a singleton}\}$.

Now, let $F^{t,\omega}\in {\cal B}(\RR^d)$ be a measurable set
with $P_t(F^{t,\omega})=1$ and such that $T^{t,\omega,n}(x)=\nabla
\varphi^{t,\omega,n}(x)$ is defined for all $x\in F^{t,\omega}$.
Then,  on $F^{t,\omega}\cap [\tilde{\Omega}_T]_{(t,\omega)}$ it
must hold that
$$T^n(t,\omega,x)= T^{t,\omega,n}(x) = \nabla
\varphi^{t,\omega,n}(x),$$ and we conclude that for all
$(t,\omega)\in \Omega'_T$,
\begin{equation*}
\begin{split}
\pi^{t,\omega,n}(C\times D)  = &  \int_{\RR^d} \1_{C\cap
F^{t,\omega}\cap
[\tilde{\Omega}_T]_{(t,\omega)}}(x)\1_D(T^{t,\omega,n}(x))P_t(dx),
\\
= & \int_{\RR^d} \1_{C\cap F^{t,\omega}\cap
[\tilde{\Omega}_T]_{(t,\omega)}}(x)\1_D(T^n(t,\omega,x))P_t(dx) \\
= & \int_{\RR^d} \1_{C}(x)\1_D(T^n(t,\omega,x))P_t(dx)
\end{split}
\end{equation*}

\Box

We point out that  Theorem \ref{theo:transportprocess} implies

\begin{coro} $T^{n}(t,\omega,x)=T^{t,\omega,n}(x)$ holds $\PP(d\omega)\otimes dt \otimes P_t(dx)\mbox{-almost surely }
$. Consequently, $T^{t,\omega,n}(x)$ is measurable with respect to
the completed  $\sigma-$field of ${\cal P}red^n\otimes {\cal
B}(\RR^d)$ with respect to $\PP(d\omega)\otimes dt \otimes
P_t(dx)$.
\end{coro}

\section{Pathwise convergence and rates for  stochastic
particle systems to Landau process}

{\bf Proof of Proposition \ref{prop:Brms}}

From the proof of Theorem \ref{theo:transportprocess}, it is clear
that integrals with respect to the measures $\1_{A^{k,n}\cap
\tilde{\Omega}_T} P_t(dx)\otimes dt$ and $\1_{A^{k,n}}
P_t(dx)\otimes dt$ are indistinguishable. By considering quadratic
variations, the same is seen to hold for the stochastic integrals
with respect to $\1_{A^{k,n}\cap \tilde{\Omega}_T} W_P^i(dx,dt)$
and $\1_{A^{k,n}} W_P^i(dx,dt)$. Write
$$ B_t^{ik,n,m}$$
for the $m-$th coordinate of the process $B_t^{ik,n}$ in
\eqref{eq:lesbrowns}, which is a real valued continuous local
martingale with respect to ${\cal F}_t^n$ (see \cite{Walsh:84}).
Then, we have that
\begin{equation*}
\begin{split}
\langle B^{ik,n,m},B^{i'k',n,m'}\rangle_t(\omega)= & n
\delta_{(i,m),(i',m')}\int_0^t \int_{\RR^d} \1_{A^{k,n}\cap
A^{k',n} \cap \tilde{\Omega}_T}(s,\omega,x) P_s(dx) ds \\
& = n\delta_{(i,k,m),(i',k',m')}\int_0^t \int_{\RR^d}
\1_{A^{k,n}\cap \tilde{\Omega}_T }(s,\omega,x) P_s(dx) ds, \\
\end{split}
\end{equation*}
by  step (a) in the proof of Theorem \ref{theo:transportprocess}.
 Now, for $(s,\omega)\in \Omega'_T$ the points
$X^1_s(\omega),\dots,X^n_s(\omega)$ are all different, and
consequently we have that
\begin{equation*}
\begin{split}
\int_{\RR^d} \1_{A^{k,n}\cap \tilde{\Omega}_T }(s,\omega,x)
P_s(dx) = & P_s(\{x:
T^n(s,\omega,x)=X^k_s(\omega)\}) \\
 = & P_s(\{x:
T^{s,\omega,n}(x)=X^k_s(\omega)\}) \\
= & \pi^{s,\omega,n} (\{(x,y): y= X^k_s(\omega)\}) \\
= & \nu_s^n(X^k_s(\omega)) \\
= & \frac{1}{n} \\
\end{split}
\end{equation*}
Thus, we have $ \langle B^{ik,n,m},B^{i'k',n,m'}\rangle_t=
t\delta_{(i,k,m),(i',k',m')},$ and the result follows.

\Box

\medskip

 We now are ready  to prove  Theorem \ref{main}.
\medskip

{\bf Proof of Theorem \ref{main}, {\it a)}} Let us fix $n\in
\NN^*$, and define for $i=1,\dots n$,
\begin{equation*}
X^{i,n}_t=X_0^i+\frac{1}{\sqrt{n}}\int_0^t \sum_{k=1}^n
\sigma(X^{i,n}_s-X^{k,n}_s) dB_s^{ik,n}+\frac{1}{n}\int_0^t
\sum_{k=1}^n b(X^{i,n}_s-X^{k,n}_s) ds
\end{equation*}
or equivalently, in an indistinguishable way,
\begin{equation*}
\begin{split}
X^{i,n}_t=&X_0^i+\int_0^t \int_{\RR^d} \sum_{k=1}^n
\sigma(X^{i,n}_s-X^{k,n}_s) \1_{A^{k,n}}(s,y)W_P^i(dy, ds)
\\ & +\int_0^t \int_{\RR^d} \sum_{k=1}^n
b(X^{i,n}_s-X^{k,n}_s)\1_{A^{k,n}}(s,y)P_s(dy) ds \\
\end{split}
\end{equation*}

By standard arguments and the fact that the sets $A^{k,n}$ are
disjoint (step (a) of the proof of Theorem
\ref{theo:transportprocess}), we have
\begin{equation}\label{EXin-Xi}
\begin{split}
E\left(|X^{i,n}_t-X^i_t|^2\right)\leq & \int_0^t \int_{\RR^d} E
\left(\sum_{k=1}^n \left(
\left[\sigma(X^{i,n}_s-X^{k,n}_s)-\sigma(X^i_s-y)\right]^2\1_{A^{k,n}}(s,y)
\right)\right)
 P_s(dy) ds \\
 & + \int_0^t \int_{\RR^d} E \left(\sum_{k=1}^n
\left(
\left[b(X^{i,n}_s-X^{k,n}_s)-b(X^i_s-y)\right]^2\1_{A^{k,n}}(s,y)
\right)\right)
 P_s(dy) ds
\end{split}
\end{equation}
The first term in the right hand side of (\ref{EXin-Xi}) is
bounded by
\begin{equation*}
\begin{split}
& C\int_0^t \int_{\RR^d} E \left(\sum_{k=1}^n \left(
\left[\sigma(X^{i,n}_s-X^{k,n}_s)-\sigma(X^i_s-X^{k,n}_s
)\right]^2\1_{A^{k,n}}(s,y) \right)\right)
 P_s(dy) ds \\
& +C \int_0^t \int_{\RR^d} E\left( \sum_{k=1}^n \left(
\left[\sigma(X^i_s-X^{k,n}_s )-
\sigma(X^i_s-T^{n}(s,y))\right]^2\1_{A^{k,n}}(s,y) \right)\right)
 P_s(dy) ds \\
& +C \int_0^t \int_{\RR^d} E \left(\sum_{k=1}^n \left( \left[
\sigma(X^i_s-T^{n}(s,y))-\sigma(X^i_s-y)\right]^2\1_{A^{k,n}}(s,y)
\right)\right)
 P_s(dy) ds \\
\leq & \quad C\int_0^t  E \left(\sum_{k=1}^n \left(
\left|X^{i,n}_s- X^i_s\right|^2\int_{\RR^d} \1_{A^{k,n}}(s,y)
P_s(dy) \right)\right)
ds \\
& +C \int_0^t \int_{\RR^d} E \left(\sum_{k=1}^n \left(
\left|X^{k,n}_s- T^{n}(s,y))\right|^2 \1_{A^{k,n}
}(s,y)\right)\right)
 P_s(dy) ds \\
 & +C \int_0^t \int_{\RR^d} E \left(\sum_{k=1}^n \left(
\left|T^{n}(s,\omega,y)-y \right|^2 \1_{A^{k,n}}(s,y)
\right)\right)
 P_s(dy) ds \\
 &=  \quad C\int_0^t  E\left(  \left|X^{i,n}_s-
X^i_s\right|^2 \right)ds \\
& + C \int_0^t  E \left(\sum_{k=1}^n \left( \left|X^{k,n}_s- X^k_s
\right|^2 \int_{\RR^d} \1_{A^{k,n}}(s,y) P_s(dy) \right)\right)
  ds \\
  & +C \int_0^t \int_{\RR^d} E\left(
\left|T^{n}(s,\omega,y)-y \right|^2\right)
 P_s(dy) ds \\
  &=  \quad C\int_0^t  E  \left(\left|X^{i,n}_s-
X^i_s\right|^2 \right)ds + C \int_0^t \frac{1}{n} E
\left(\sum_{k=1}^n
\left|X^{k,n}_s- X^k_s \right|^2\right) ds \\
 & +C \int_0^t \int_{\RR^d} E\left(
\left|T^{n}(s,\omega,y)-y \right|^2\right)
 P_s(dy) ds \\
& = 2 C\int_0^t  E \left(\left|X^{i,n}_s- X^i_s\right|^2\right) ds
+ C \int_0^t E \left(W_2^2(\nu_s^n,P_s)\right)ds \\\end{split}
\end{equation*}
by exchangeability of $((X^{1,n},X^1),\dots,(X^{n,n},X^n))$.  A
similar bound is obtained for the second term in (\ref{EXin-Xi}).
We deduce by Gronwall's lemma that
\begin{equation*}
E\left(|X^{i,n}_t-X^i_t|^2\right)\leq C\exp(C'T) \int_0^t E
(W_2^2\left(\nu_s^n,P_s)\right)ds
\end{equation*}
By a little finer argument using a Burkholder-Davis-Gundy
inequality, we can obtain as usual an estimate of the form
\begin{equation*}
E\left(\sup_{t\in [0,T]}|X^{i,n}_t-X^i_t|^2\right)\leq C\exp(C'T)
\int_0^T E \left(W_2^2(\nu_s^n,P_s)\right)ds
\end{equation*}

\Box

\bigskip

 We recall a result
proved in Rachev and R\"uschendorf \cite{Rachev:98} giving
$L^2$-rates of convergence of empirical measures in the
Wasserstein metric.

\begin{theo} (\cite{Rachev:98} Theorem 10.2.1)
\label{wass-empirmeas} Let $\mu$ a probability on $\RR^d$ and let
$Y^1, Y^2, \ldots, Y^n$ be independent identically distributed
random variables with law $\mu$. Let $\mu_n$ be the empirical
measure of these variables. Then, if $\mu$ has high enough finite
absolute moments: $c:=\int_{\RR^d} |y|^{d+5} \mu(dy) <\infty,$
there is a constant $C$ depending only on $c$ and on the dimension
$d$, such that
$$E\left( W_2^2(\mu_n,\mu)\right)\leq C n^{{-2\over d+4}}.$$
\end{theo}

Denote by ${\cal W}_2$ the Wasserstein distance between
probability measures $Q$ on the path space
$\mathcal{C}_T:=C([0,T],\RR^d)$, such that $\int_{\mathcal{C}_T}
\sup_{0\leq t \leq T} |x(t)|^2Q(dx)<\infty$.

From the previous result and Lemma \ref{lem:momboud}, it is simple
to deduce the following
\begin{coro}\label{coro:weakrate}
Let $P$ be the pathwise law of the nonlinear process
\eqref{eq:nonlinproc}. Under the  assumptions of Theorem
\ref{main} and  moreover that $\int_{\RR^d} |y|^{d+5}
P_0(dy)<\infty$, we have that
\begin{equation*}
{\cal W}_2^2(law(X^{1,n}),P)\leq C_{T,d}n^{{-2\over d+4}}.
\end{equation*}
\end{coro}

\Box

\medskip

The previous results are the first convergence rates obtained  so
far for stochastic particle systems of the ``Landau type''
\eqref{eq:systparticules}, and they are not specific to the
particular coefficients of the Landau equation
\eqref{eq:Landaueq}. They justify the interest of the particle
systems introduced in \eqref{eq:systparticules} and are the first
step in the construction and the numerical study of a simulation
algorithm for $(P_t)_t$. We notice that since we deal with
space-time random fields, the dependence of the results on the
dimension $d$ is somewhat expectable, as opposite to the situation
in the McKean-Vlasov model. The techniques we have introduced
provide some insight about that dependence.

\bigskip

{\bf Acknowledgements} The authors are very grateful to Roberto
Cominetti for helpful suggestions about the theory of set-valued
mappings.

\end{document}